\documentclass[11pt]{article}

\usepackage{arxiv}

\usepackage[utf8]{inputenc} 
\usepackage[T1]{fontenc}    
\usepackage{lmodern}

\usepackage{amsmath,amssymb,amsthm}
\usepackage{amsfonts}       
\usepackage{nicefrac}       

\usepackage{graphicx}
\usepackage{booktabs}       
\usepackage{float}
\usepackage{xcolor}
\usepackage{microtype}      
\usepackage{lipsum}         

\usepackage[hyphens]{url}   
\usepackage{cite}           

\usepackage{hyperref}

\usepackage{orcidlink}      
\usepackage{doi}

\title{Combinatorial and Sampling Approaches to Determining the Priority Vector from an Incomplete Pairwise Comparison Matrix on the Set of Spanning Trees}

\date{July 1, 2026}	

\author{
{ \orcidlink{0000-0002-0821-4877} Vitaliy Tsyganok} \\
\small Institute for Information Recording of the National Academy of Sciences of Ukraine, Kyiv, Ukraine\\
\small Taras Shevchenko National University of Kyiv, Kyiv, Ukraine\\
\small National Technical University of Ukraine ``Igor Sikorsky Kyiv Polytechnic Institute'', Kyiv, Ukraine\\
\texttt{tsyganok@ipri.kiev.ua} \\
\And
{ \orcidlink{0000-0002-6117-5846} Oksana Mulesa} \\
\small University of Pre\v{s}ov, Pre\v{s}ov, Slovakia\\
\small Uzhhorod National University, Uzhhorod, Ukraine\\
\texttt{oksana.mulesa@unipo.sk} \\
}

\renewcommand{\shorttitle}{\textit{Tsyganok, Mulesa} \small Combinatorial and Sampling Approaches for Set of Spanning Trees}

\hypersetup{
pdftitle={Combinatorial and Sampling Approaches for Incomplete PCMs on Spanning Trees},
pdfsubject={SC.AI},
pdfauthor={Vitaliy~Tsyganok, Oksana~Mulesa},
pdfkeywords={pairwise comparison matrix, priority vector, spanning tree, Monte Carlo, combinatorial method},
}

\begin{document}
\maketitle

\begin{abstract}
This paper considers combinatorial and sampling approaches to deriving a priority vector from an incomplete pairwise comparison matrix. It is shown that the informationally significant combinations of matrix elements are in one-to-one correspondence with spanning trees of the graph induced by the matrix. The combinatorial approach is based on the enumeration of all spanning trees and the componentwise aggregation of the priority vectors derived from them. To reduce computational complexity, a sampling Monte Carlo approach with a sequential stopping criterion is proposed, ensuring a prescribed accuracy for estimating the components of the resulting vector.
\end{abstract}

\noindent\textbf{Keywords:} pairwise comparison matrix, priority vector, spanning tree, Monte Carlo, combinatorial method.

\section{Introduction}\label{sect1}

The combinatorial approach to deriving a priority vector from a pairwise comparison matrix (PCM) was first proposed in~\cite{TsyganokComb}. It is based on the complete enumeration of all informationally significant combinations of matrix elements in order to recover a consistent comparison structure. The main idea is that not all matrix elements are independent: if the matrix is ideally consistent, then a part of its entries can be uniquely reconstructed from the others.

For an ideally consistent PCM \(M=\{m_{ij}\}\), there exists a set of alternative scores \(\langle A_1, A_2, \dots, A_n \rangle\) such that each matrix element is determined by a function of the corresponding alternative scores:
\[
m_{ij}=f(A_i,A_j).
\]
Depending on the adopted scale, one may use:
\begin{itemize}
    \item multiplicative comparisons
    \[
    m_{ij}=\frac{A_i}{A_j},
    \]
    \item additive comparisons
    \[
    m_{ij}=A_i-A_j.
    \]
\end{itemize}

A key consequence is the consistency relation between matrix elements. In particular, for any \(i,j,k\):
\begin{itemize}
    \item in the multiplicative case,
    \[
    m_{ij}=\frac{m_{ik}}{m_{jk}},
    \]
    \item in the additive case,
    \[
    m_{ij}=m_{ik}-m_{jk}.
    \]
\end{itemize}

This relation makes the combinatorial approach possible: instead of processing the whole matrix at once, one may identify minimal subsets of elements containing sufficient information to reconstruct all remaining entries. In graph-theoretic terms, alternatives correspond to vertices and available pairwise comparisons correspond to edges. Then the problem reduces to finding minimal connected structures that allow one to reconstruct all missing comparisons uniquely. The minimum number of informationally significant elements is \(n-1\), which corresponds to a spanning tree of the associated graph.

Hence, the essence of the combinatorial method is the enumeration of all informationally significant subsets of elements, i.e., all spanning trees of the graph corresponding to the PCM, the construction of a consistent matrix or, equivalently, a uniquely defined priority vector on each such structure, and the subsequent aggregation of these vectors into a generalized priority vector. This approach yields an effective solution~\cite{TsyganokEff2010}, but it requires substantial computational resources due to the combinatorial growth in the number of admissible subsets.

Subsequent studies established the mathematical equivalence between the results obtained by the combinatorial method and those obtained via the row geometric mean of a complete multiplicative PCM~\cite{MatEquiv}. More general results for incomplete multiplicative and additive PCMs were obtained in~\cite{BozokiTsyganok}, where the optimality of the combinatorial method was shown in the framework of (logarithmic) least squares. In the same line of research, the computational complexity of the combinatorial method was compared with that of an equivalent optimization-based approach~\cite{BozokiMethod}.

If the graph is complete, one has to enumerate all spanning trees of the complete graph; if the PCM is incomplete, only the spanning trees contained in the corresponding subgraph of available comparisons are considered. This setting is particularly important in practical applications, where some pairwise comparisons may be unavailable due to time limitations, elicitation cost, or the impossibility of obtaining certain expert judgments.

The main drawback of the combinatorial method is that the number of spanning trees grows very rapidly with problem size, so a complete enumeration of all structures requires substantial computational time and memory. In addition, for each spanning tree one has not only to generate the structure itself, but also to derive the corresponding priority vector and then aggregate all such vectors. For incomplete graphs, one must furthermore verify that the generated tree belongs to the original graph. These factors make exhaustive combinatorial enumeration impractical for medium- and large-scale problems.

This motivates the transition to a sampling approach based on the random generation of spanning trees and the statistical estimation of the resulting priority vector~\cite{MetropolisUlam1949,Rubinstein}. Let
\[
T^{(1)},T^{(2)},\dots,T^{(N)}
\]
be randomly selected spanning trees of a graph \(G\). Instead of exact aggregation over the whole set of spanning trees, one can use the sample estimator
\[
\hat w^{(N)}=\frac{1}{N}\sum_{k=1}^{N} w\!\left(T^{(k)}\right).
\]
For sufficiently large \(N\) and under a proper sampling scheme, this estimator approximates the exact combinatorial result. Moreover, if a sequential sampling procedure is employed and the estimation accuracy of each component is monitored, the computation can be terminated as soon as the required precision is achieved~\cite{Wald1945}.

Thus, this paper considers a Monte Carlo sampling approach to the componentwise aggregation of priority vectors over the set of spanning trees of the graph associated with a PCM. This approach preserves the conceptual connection with the combinatorial method while significantly reducing computational cost by replacing exhaustive enumeration with random sampling equipped with a sequential stopping rule.

\section{Problem statement}\label{sect2}

Consider an incomplete pairwise comparison matrix associated with a connected undirected graph
\[
G=(V,E),
\]
where the vertex set
\[
V=\{1,2,\dots,n\}
\]
corresponds to the alternatives, and the edge set \(E\) corresponds to the available pairwise comparisons. The graph \(G\) is assumed to be connected, since only in this case does a spanning tree exist and, consequently, a consistent relation among all alternatives can be reconstructed.

Let
\[
\mathcal{T}(G)=\{T_1,T_2,\dots,T_M\}
\]
denote the set of all spanning trees of \(G\), where
\[
M=|\mathcal{T}(G)|.
\]

Each spanning tree \(T\in\mathcal{T}(G)\) determines a priority vector
\[
w(T)=\bigl(w_1(T),w_2(T),\dots,w_n(T)\bigr),
\]
subject to
\[
w_i(T)>0,\qquad i=1,\dots,n,
\]
\[
\sum_{i=1}^{n}w_i(T)=1.
\]

In the multiplicative model, if the edge \((i,j)\in T\) is assigned the pairwise comparison value
\[
a_{ij}=\frac{v_i}{v_j},
\]
then all unnormalized scores \(v_i\) can be reconstructed from the tree \(T\) up to a common positive factor. Normalization yields
\[
w_i(T)=\frac{v_i}{\sum_{j=1}^{n}v_j},\qquad i=1,\dots,n.
\]

Thus, each spanning tree defines a minimally sufficient information structure for constructing one consistent priority vector. The generalized combinatorial result is defined as the componentwise mean of all such vectors:
\[
\bar w=\frac{1}{M}\sum_{T\in\mathcal{T}(G)}w(T).
\]

Componentwise, this means
\[
\bar w_i=\frac{1}{M}\sum_{T\in\mathcal{T}(G)}w_i(T),\qquad i=1,\dots,n.
\]

The goal is to construct a method that estimates the vector \(\bar w\) without full enumeration of all spanning trees, using a random sample of trees together with a sequential stopping criterion that guarantees a prescribed estimation accuracy for each component of the resulting vector.

\section{Combinatorial approach}\label{sect3}

The combinatorial approach consists of the explicit enumeration of all spanning trees of \(G\), the construction of a priority vector for each tree, and the subsequent componentwise aggregation of the obtained vectors.

\subsection{Enumeration of spanning trees}

If the graph \(G\) is complete, i.e., \(G=K_n\), then the number of spanning trees is
\[
n^{n-2}
\]
by Cayley's formula~\cite{Cayley}. For the systematic generation of these trees, Prüfer codes are convenient~\cite{Prufer}, since each code of length \(n-2\) corresponds uniquely to a spanning tree of the complete graph.

In the incomplete-graph case, Prüfer codes may still be used to generate candidate trees, followed by a verification of whether all generated edges belong to the set \(E\). Alternatively, one may employ specialized algorithms for generating spanning trees directly in the given graph.

\subsection{Construction of a priority vector on a spanning tree}

For each spanning tree \(T\in\mathcal{T}(G)\), a consistent priority vector \(w(T)\) is constructed from the pairwise comparison values associated with its edges. Since a spanning tree is a connected acyclic structure, it contains exactly \(n-1\) independent relations among the alternatives, which is sufficient to recover all relative scores up to a common multiplicative factor.

After reconstructing the unnormalized scores \(v_1,\dots,v_n\), one performs normalization:
\[
w_i(T)=\frac{v_i}{\sum_{j=1}^{n}v_j},\qquad i=1,\dots,n.
\]

The resulting vector \(w(T)\) is interpreted as the priority vector induced by the informationally significant combination of pairwise comparisons corresponding to the spanning tree \(T\).

\subsection{Aggregation over the set of spanning trees}

After constructing priority vectors for all spanning trees, the combinatorial result is defined as their componentwise mean:
\[
\bar w=\frac{1}{|\mathcal{T}(G)|}\sum_{T\in\mathcal{T}(G)}w(T).
\]

This approach is theoretically complete because it accounts for the entire set of minimally sufficient information structures. However, it has a substantial computational cost. If each spanning tree requires the reconstruction of a priority vector and an update of the aggregated result, the total effort grows proportionally to the number of spanning trees, which becomes very large for graphs of medium and large size.

\subsection{Advantages and limitations}

The main advantages of the combinatorial approach are:
\begin{itemize}
    \item completeness with respect to the whole set of minimally informationally sufficient structures;
    \item a natural graph-theoretic interpretation through spanning trees of the PCM graph;
    \item uniqueness of the aggregated result once the aggregation rule is fixed.
\end{itemize}

Its main limitations are:
\begin{itemize}
    \item the excessive number of spanning trees as \(n\) grows;
    \item the computational cost of generating trees and constructing priority vectors;
    \item the need to test admissibility of trees in the incomplete-graph case;
    \item the practical infeasibility of exhaustive enumeration for large-scale problems.
\end{itemize}

These limitations motivate the transition to a sampling-based approach.

\section{Sampling approach}\label{sect4}

The sampling approach is based on the random generation of spanning trees from \(\mathcal{T}(G)\) and on the estimation of the aggregated priority vector from a sample of such trees.

Let
\[
T^{(1)},T^{(2)},\dots,T^{(N)}
\]
be randomly generated spanning trees of the graph \(G\). For each tree \(T^{(k)}\), a priority vector is computed:
\[
w\!\left(T^{(k)}\right)=\bigl(w_1(T^{(k)}),\dots,w_n(T^{(k)})\bigr).
\]

Then the sample estimator of the resulting vector is
\[
\hat w^{(N)}=\frac{1}{N}\sum_{k=1}^{N}w\!\left(T^{(k)}\right),
\]
that is, componentwise,
\[
\hat w_i^{(N)}=\frac{1}{N}\sum_{k=1}^{N}w_i(T^{(k)}),\qquad i=1,\dots,n.
\]

If the spanning trees are generated independently and uniformly over \(\mathcal{T}(G)\), then each component \(\hat w_i^{(N)}\) is an unbiased estimator:
\[
\mathbb{E}\bigl[\hat w_i^{(N)}\bigr]=\bar w_i.
\]

Therefore, the sampling approach can be regarded as a statistically valid substitute for full combinatorial enumeration in problems where the exact computation of \(\bar w\) is too expensive.

\subsection{Stopping criterion}

Unlike fixed-sample procedures, the sequential sampling approach determines the number of generated trees adaptively. Let, for each component \(i\), the sample variance after \(N\) observations be
\[
S_{i,N}^2=
\frac{1}{N-1}
\sum_{k=1}^{N}
\left(
w_i(T^{(k)})-\hat w_i^{(N)}
\right)^2.
\]

Then the half-width of the confidence interval for the \(i\)-th component is
\[
\Delta_{i,N}=t_{1-\alpha/2,\;N-1}\frac{S_{i,N}}{\sqrt N},
\]
where \(t_{1-\alpha/2,\;N-1}\) is the corresponding Student quantile and \(1-\alpha\) is the confidence level.

If \(\varepsilon>0\) is a prescribed absolute error tolerance, then a natural stopping condition is
\[
\Delta_{i,N}\le \varepsilon,\qquad i=1,\dots,n.
\]

Equivalently,
\[
\max_{1\le i\le n}\Delta_{i,N}\le \varepsilon.
\]

Thus, the generation of spanning trees continues until all components of the resulting vector are estimated with the required accuracy. This type of adaptive sample size control follows the general idea of sequential statistical estimation~\cite{Wald1945}.

\subsection{Assumptions on tree generation}

A key condition for the correctness of the sampling approach is the tree generation scheme. If the graph is complete, uniform generation can be implemented via random Prüfer codes. In the incomplete-graph case, the problem is more difficult because not every Prüfer code yields an admissible tree. In this situation one may use:
\begin{itemize}
    \item generation of candidate trees followed by rejection of inadmissible ones;
    \item specialized algorithms for uniform generation of spanning trees in a given graph;
    \item weighted estimation schemes in the case of non-uniform sampling.
\end{itemize}

In the sequel, it is assumed that the generation of spanning trees is uniform or, more generally, that the adopted sampling scheme provides statistically correct estimation of the components of \(\bar w\).

\section{Sequential Monte Carlo algorithm}\label{sect5}

Below we present a general scheme of a sequential Monte Carlo algorithm for estimating the aggregated priority vector over the set of spanning trees.

\subsection{Input data}

The algorithm requires:
\begin{itemize}
    \item a connected graph \(G=(V,E)\) corresponding to an incomplete PCM;
    \item a procedure for random generation of a spanning tree \(T\in\mathcal{T}(G)\);
    \item a procedure for constructing the priority vector \(w(T)\) from a spanning tree \(T\);
    \item the significance level \(\alpha\);
    \item the absolute error tolerance \(\varepsilon\);
    \item an initial sample size \(N_0\ge 2\).
\end{itemize}

\subsection{Recursive estimation}

To reduce memory usage, it is convenient not to store all constructed priority vectors but to update the componentwise statistics recursively. Suppose that after \(N-1\) steps, for each component \(i\), one knows:
\begin{itemize}
    \item the current mean \(m_{i,N-1}\);
    \item the auxiliary quantity \(M2_{i,N-1}\) used for variance computation.
\end{itemize}

After receiving a new value \(w_i^{(N)}=w_i(T^{(N)})\), perform:
\[
\delta_{i,N}=w_i^{(N)}-m_{i,N-1},
\]
\[
m_{i,N}=m_{i,N-1}+\frac{\delta_{i,N}}{N},
\]
\[
M2_{i,N}=M2_{i,N-1}+\delta_{i,N}\bigl(w_i^{(N)}-m_{i,N}\bigr).
\]

Then
\[
\hat w_i^{(N)}=m_{i,N},
\qquad
S_{i,N}^2=\frac{M2_{i,N}}{N-1}.
\]

\subsection{Algorithm description}

The sequential algorithm can be described as follows:
\begin{enumerate}
    \item Set the parameters \(\alpha\), \(\varepsilon\), and \(N_0\).
    \item Generate \(N_0\) random spanning trees.
    \item Construct a priority vector for each generated tree.
    \item Update the componentwise sample means and variances.
    \item For each component compute
    \[
    \Delta_{i,N}=t_{1-\alpha/2,\;N-1}\frac{S_{i,N}}{\sqrt N}.
    \]
    \item If, for all \(i=1,\dots,n\),
    \[
    \Delta_{i,N}\le \varepsilon,
    \]
    then terminate the computation.
    \item Otherwise, generate a new spanning tree, compute the corresponding priority vector, update the statistics, and return to Step~5.
\end{enumerate}

\subsection{Properties of the algorithm}

The proposed algorithm has the following properties:
\begin{itemize}
    \item the resulting estimate of the priority vector is the componentwise sample mean;
    \item the sample size is determined adaptively according to the actual variability of the components;
    \item the algorithm does not require storing the whole sample;
    \item the stopping rule is directly linked to the prescribed estimation accuracy.
\end{itemize}

Compared with full enumeration of all spanning trees, the proposed approach can substantially reduce computational effort, especially in large-scale problems.

\section{Conclusions}\label{sect6}

This paper considered the problem of deriving a priority vector from an incomplete pairwise comparison matrix over the set of spanning trees of the corresponding graph. It was shown that the minimally informationally significant combinations of PCM elements naturally correspond to spanning trees of the graph induced by the matrix. For each such tree, a consistent priority vector can be constructed, and the generalized result is defined as the componentwise aggregation of the vectors obtained over the whole set of spanning trees.

The combinatorial approach was analyzed as a complete method based on exhaustive enumeration of all minimally sufficient information structures, but with high computational complexity. As an alternative, a Monte Carlo sampling approach was proposed, in which the resulting priority vector is estimated from a random sample of spanning trees. For this approach, a componentwise sequential stopping criterion was formulated, allowing one to terminate the computation once the prescribed accuracy has been achieved for all components.

The proposed formulation provides a basis for the practical use of the sampling approach in problems where full combinatorial enumeration is computationally prohibitive. Further research should focus on efficient algorithms for uniform generation of spanning trees in incomplete graphs, as well as on experimental comparisons of the accuracy and computational efficiency of the combinatorial and sampling approaches.


\begin{thebibliography}{99}

\bibitem{TsyganokComb}
Tsyganok, V. 2000.
Combinatorial Method of Pairwise Comparisons with Feedback.
\emph{Data Recording, Storage \& Processing} 2, 92--102 (in Ukrainian).
Available at: \url{https://dss-lab.org.ua/documents/Tsyganok_DRSP-2000-2(ukr).pdf}.

\bibitem{TsyganokEff2010}
Tsyganok, V. 2010.
Investigation of the Aggregation Effectiveness of Expert Estimates Obtained by the Pairwise Comparison Method.
\emph{Mathematical and Computer Modelling} 52 (3--4), 538--544.
DOI: \url{https://doi.org/10.1016/j.mcm.2010.03.010}.

\bibitem{MatEquiv}
Lundy, M., Siraj, S., and Greco, S. 2017.
The Mathematical Equivalence of the ``Spanning Tree'' and Row Geometric Mean Preference Vectors and its Implications for Preference Analysis.
\emph{European Journal of Operational Research} 257 (1), 197--208.
DOI: \url{https://doi.org/10.1016/j.ejor.2016.07.029}.

\bibitem{BozokiTsyganok}
Boz\'oki, S., and Tsyganok, V. 2019.
The (logarithmic) least squares optimality of the arithmetic (geometric) mean of weight vectors calculated from all spanning trees for incomplete additive (multiplicative) pairwise comparison matrices.
\emph{International Journal of General Systems} 48 (4), 362--381.
DOI: \url{https://doi.org/10.1080/03081079.2019.1585432}.

\bibitem{BozokiMethod}
Boz\'oki, S., F\"ul\"op, J., and R\'onyai, L. 2010.
On Optimal Completion of Incomplete Pairwise Comparison Matrices.
\emph{Mathematical and Computer Modelling} 52 (1--2), 318--333.
DOI: \url{https://doi.org/10.1016/j.mcm.2010.02.047}.

\bibitem{MetropolisUlam1949}
Metropolis, N., and Ulam, S. 1949.
The Monte Carlo method.
\emph{Journal of the American Statistical Association} 44 (247), 335--341.
DOI: \url{https://doi.org/10.1080/01621459.1949.10483310}.

\bibitem{Rubinstein}
Rubinstein, R.Y., and Kroese, D.P. 2017.
\emph{Simulation and the Monte Carlo Method}.
3rd ed.
Hoboken, NJ: Wiley.

\bibitem{Wald1945}
Wald, A. 1945.
Sequential tests of statistical hypotheses.
\emph{The Annals of Mathematical Statistics} 16 (2), 117--186.
DOI: \url{https://doi.org/10.1214/aoms/1177731118}.

\bibitem{Cayley}
Cayley, A. 1889.
A theorem on trees.
\emph{Quarterly Journal of Mathematics} 23, 376--378.

\bibitem{Prufer}
Pr\"ufer, H. 1918.
Neuer Beweis eines Satzes \"uber Permutationen.
\emph{Archiv der Mathematik und Physik} 27, 742--744.

\end{thebibliography}
\end{document}